\documentclass{amsart}
\usepackage[cp1251]{inputenc}
\usepackage[english,russian]{babel}
\usepackage{amsmath}
\usepackage{amssymb}
\usepackage{amsfonts}
\def\udcs{517.518.5} %Здесь автор определяет УДК своей работы

\newtheorem{theorem}{Теорема}

\begin{document}
УДК \udcs
\thispagestyle{empty}

\title[Об оценке тригонометрических интегралов с квадратичной фазой]
{Об оценке тригонометрических интегралов с квадратичной фазой}
% Указываем название статьи.
%Сокращенное название указывается в квадратных скобках для колонтитулов,
%если полное название не помещается в строку

\author{ И.~А.~Икромов, А.Р.Сафаров, А.Т.Абсаламов}
%Указываем авторов
\address{Исроил Акрамович Икромов, %Имя, Отчество, Фамилия первого автора
\newline\hphantom{iii} Акбар Рахманович Сафаров,
% Имя, Отчество, Фамилия второго автора
\newline\hphantom{iii} Акмал Толлибоевич Абсаламов,
\newline\hphantom{iii} Институт математики имени В.И.Романовского, 
\newline\hphantom{iii} Самаркандский государственный университет % Место работы
\newline\hphantom{iii} ул.Университетский бульвар, 15, % Адрес (улица, дом, строение и т.п.)
\newline\hphantom{iii} 140104, г. Самарканд, Узбекистан}%  Адрес (почтовый индекс, город, страна)
\email{safarov-akbar@mail.ru}% Ваш электронный адрес для переписки

\thanks{\sc Ikromov I.A., Safarov A.R., Absalamov A.T. %  Ф.И.О. авторов на английском языке
On estimates for trigonometric integrals with quadratic phase}% название статьи на английском языке
\thanks{\copyright \ 2021 Икромов И.А., Сафаров А.Р., Абсаламов А.Т. }
%\thanks{\rm Работа поддержана РФФИ (грант 03-01-11111)}

% (указываем дату отправки, строка будет при получении изменена)

%%%%%%%%%%%%%%%%%%%%%%%%%%%%%%%%%%%%%%%%%%%%%%%%%%%%%%%%%%%%%%%%
\maketitle
{
\small
\begin{quote}
\noindent{\bf Аннотация. } В статье рассматривается проблема суммируемости для тригонометрических интегралов с квадратичной фазой. Аналогичная задача рассмотрена в работах \cite{Chub}, \cite{Chax}, \cite{Jabbar} в частных случаях.
Наши результаты обобщают результаты этих работ на кратные тригонометрические интегралы.
\medskip

\noindent{\bf Abstract. } In paper this paper it is considered the summation problem for trigonometric
integrals with quadratic phase. This problem considered in the
papers \cite{Chub},\cite{Chax},\cite{Jabbar} in particular cases.
Our results generalized the results of that papers to multidimensional trigonometrical integrals.
\medskip

\noindent{\bf Ключевые слова:}{ тригонометрический интеграл, экспонент, сумма, фаза, многочлен.}
\medskip

\end{quote}
}

\section{Введение}
Пусть
$P(x,s)\in{\mathbb{R}[x]}$ многочлен от $x\in{\mathbb{R}^k}$ с коэффициентами
$s\in{\mathbb{R}^N}$. Через $Q$ обозначается компактное множество в $\mathbb{R}^k.$

Рассмотрим тригонометрический интеграл
\begin{equation}\label{T(s) formulasi}
T(s)=\displaystyle\int\limits_{Q}{\exp(iP(x,s))dx}.
\end{equation}
\textbf{Постановка задачи:} Найти точную нижнюю грань $p_{0}$ чисел $p$ таких, что $T\in \mathbb{L}_{p}(\mathbb{R}^{N}).$

Эта задача впервые была рассмотрена И.М.Виноградовым \cite{V-1} в связи с проблемой аналитической теории чисел и получена
оценка сверху для $p_{0}$ в случае $k=1$. Позднее, оценка И.М.Виноградова была улучшена в работе \cite{H-L}. В работе \cite{AKC}
указано точное значение $p_{0}$ в случае $k=1$ и доказана конечность этого числа в многомерных случаях. В работе \cite{I-2}
рассмотрены оценки снизу для числа $p_{0}$ и указано его точное значение когда коэффициенты многочлена меняются в некотором
подпространстве пространства $\mathbb{R}^{N}.$ Аналогичные задачи рассмотрены в работах \cite{SafAnal}, \cite{SafAUzmat}.

Аналог этой проблемы рассмотрен в работе \cite{M}, для случая когда $Q$ есть единичный шар с центром в начале координат и $P(x,s)$ квадратичный полином удовлетворяющий некоторому условию трансверсальности.

В работе \cite{Chub} получена оценка снизу для $p_{0}$ в случае $k=2.$

В работах \cite{Chub} и \cite{Jabbar} рассмотрена аналогичная задача в случае $k=2.$ Более того в \cite{Chub}, показано, что если $P$
однородный квадратичный полином и $k=2,$ то $p_{0}=4$ в случае когда $Q=[0,1]^{2}$ точнее при $p>4$ тригонометрический интеграл сходится и при $p\leq4$ расходится.

В данной работе мы рассмотрим задачу суммируемости тригонометрических интегралов когда $k\geq1$ и получим точный показатель сходимости $p_{0}$ в случае когда $Q=[0,1]^{k}.$

В случае когда $P$ однородный многочлен степени два получим точное значение  $p_{0}$.

Пусть полином $P$ имеет вид:
$$P(x,A,b)=(Ax,x)+(b,x),$$
где ${A=(a_{lm})^{k}}_{l,m=1}$ вещественная симметричная $k\times k$ матрица, $b:=(b_1,b_2,...,b_k)\in{\mathbb{R}^k}$ и
$({\cdot},{\cdot})$ скалярное произведение векторов. Рассмотрим тригонометрический интеграл
$$T(A,b)=\displaystyle\int\limits_{\mathbb{R}^{k}}{\exp{(iP(x,A,b))\chi_{K}(x)}dx},$$
где $K-$компактное множество и $\chi_{K}(x)-$характеристическая функция множества $K.$

Рассмотрим несобственный интеграл
$$\theta=\displaystyle\int\limits_{{\mathbb{R}}^N}
{{|T(A,b)|^{p}}dbda},$$ где $db=db_1db_2...db_k$
и $da= {\prod\limits_{1{\leq{l}}{\leq{m}}{\leq{k}}}da_{lm}}.$\\
Справедлива следующая:
\begin{theorem}\label{main Th} Пусть $K$ компактное множество, тогда интеграл  $\theta$ сходится при
${p>2k+2}$ и причем если $K$ содержит внутренную точку $x^{0}$ и существует прямая $l$ проходящая через точку $x^{0}$ такая, что множество $\{l\cap K\}$ содержит лишь конечное число точек, то  при ${p\leq{2k+2}}$ интеграл расходится. Таким образом, если $K$ компактное множество с непустой внутренностью, то
${p_0=2k+2}.$
\end{theorem}
{\bf Доказательство теоремы 1.} Оценка сверху для $p_{0}$ непосредственно следует из теоремы 1 работы \cite{LeeBak}.
Рассмотрим следующее подмножество
${\Omega{(a_{11})}}$ пространства $\mathbb{R}^{N-1}$:
$${|a_{12}|+|a_{13}|+...+|a_{1k}|}<c_{1}a_{11}, \quad
-{\frac{1}{2}}<{\frac{b_{1}}{a_{11}}}<-{\frac{1}{4}}, \quad
|{a_{lj}}-{\frac{a_{1l}a_{1j}}{a_{11}}}|{\leq{c_{2}}}, \quad
|{b_{l}}-{\frac{{2}{b}_{l}{a_{1l}}}{a_{11}}}|{\leq{c_{2}}},$$ где
$l=2,...,n$ и $c_{1}$, $c_{2}$ достаточно малые фиксированные положительные числа.

\textbf{Лемма 1.}\label{Omega}
Существует положительное число $c$ такое, что для меры Лебега множества ${\Omega{(a_{11})}}$ справедливо следующее равенство:
$${\mu{({\Omega{(a_{11})}})}} =c\cdot{a^{k}_{11}}.$$

\textbf{Доказательство леммы 1.} Берем следующее отображение
$$\xi_{1l}(A,b_{1},...,b_{k})=a_{1l},$$
$$\xi^{1}(A,b_{1},...,b_{k})=b_{1},$$
$$\xi^{l}(A,b_{1},...,b_{k})=b_{l}-
\frac{2b_{1}a_{1l}}{a_{11}},$$ $$\xi_{lj}(A,b_{1}
,..., b_{k})=a_{lj}-\frac{a_{1l}a_{1j}}{a_{11}} ,$$
$${j\leq l=2,3,...,k}.$$
Оно отображает множество ${\Omega{(a_{11})}}$ на множество ${\Omega{(\xi)}}$
и Якобиан этого отбражения равен
единице. Следовательно,
$$\mu({\Omega{(a_{11})}})=\mu({\Omega{(\xi)}}).$$
Легко показать, что для множества $\Omega{(\xi)}$:
$$\mid{\xi_{12}}\mid+\mid{\xi_{13}}\mid+...+
\mid{\xi_{1k}}\mid<c_{1} \cdot{a_{11}},$$ $$-{\frac{1}{2}}<
{\frac{\xi^{1}}{a_{11}}< -{\frac{1}{4}}},$$ $${\mid{\xi^{l}}
\mid}\leq{c_2}, \quad {\mid{\xi_{lj}}\mid}\leq{c_2}, \quad
j<l=2,3,...,k,$$ получим
$${\mu{({\Omega{(\xi)}})}}=c\cdot{a^{k}_{11}}.$$
Следовательно,
$${\mu{({\Omega{(a_{11})}})}}=c\cdot{a^{k}_{11}}.$$
\textbf{Лемма 2.}\label{T(A,b)}
Существует положительное число $L$ такое, что для любого ${a_{11}}>L$ и
$(A,b){\in{\Omega{(a_{11})}}}$ для интеграла $T(A,b)$ справедливо следующее асимптотическое равенство
$$T(A,b)={\frac{c(A,b)}{a^{\frac{1}{2}}_{11}}+O\left({\frac{1}{a_{11}}}\right)}$$
при ${a_{11}}{\rightarrow}{+{\infty}},$
причем существует положительное число ${\delta}$ такое, что для любого
$(A,b){\in{{\Omega{(a_{11})}}}}$, выполняется неравенство
$$|c(A,b)|>{\delta}>0.$$

\textbf{Лемма 2 доказывается} обычным методом стационарной фазы. Отметим, что для достаточно малых $c_{1}$, $c_{2}$ при $(A,b)\in{\Omega{(a_{11})}}$ и для достаточных больших $L$, фаза имеет осцилляции только в направлении $x_{1}$
по этому, при фиксированных значениях ${x_{2},...,x_{n}{\in{[0,1]}}}$, невырожденная критическая точка
${x_{1}}(A,b,x_{2},...,x_{n})$ лежит внутри $(0,1).$

Наконец, для интеграла ${\theta}$, имеем оценку снизу
$${\theta}{\geq}{\displaystyle\int\limits_{L}^{\infty}}
{\displaystyle\int\limits_{{\Omega{(a_{11})}}}{{|T(A,b)|^{p}}dbda}}
{\geq}{\delta}{c}{\displaystyle\int\limits_{L}^{\infty}}{a^{k-
{\frac{p}{2}}}_{11}}da_{11}.$$

Таким образом при ${p\leq{2k+2}}$ последний интеграл расходится. Теорема 1 доказана.

\section{Случай, когда $P$ однородный многочлен второй степени}

Теперь предположим, что $P(x,A)=(Ax,x)$. В работе \cite{Jabbar} доказано, что если $Q$ квадратичный полином в $\mathbb{R}^{2},$ то при $p>4$ интеграл $\theta$ сходится и при $p_{0}\leq4$ интеграл $\theta$ расходится. В данной работе мы распространяем результаты И.Ш.Джаббарова на случай, когда $Q$ многогранник в $\mathbb{R}^{2k}.$

Под многогранником мы подразумеваем конечное объединение невырожденных симплексов \cite{LebedevDisser}.
\begin{theorem}\label{Th2}
Если $P(x,A)=(Ax,x)$ и $Q$ многогранник, то при $p>2k$ интеграл $\theta$ сходится.
Если $Q=[0,1]^{k},$ то при $p\leq2k$ интеграл $\theta$ расходится.
\end{theorem}

\textbf{ Замечание 1.} В этом случае мы не можем применить результаты работы \cite{LeeBak}, так как соответствующее множество $\{x_{i}x_{j}\}_{i\leq j=1}^{n}$ не является гладкой поверхностью.

\textbf{ Замечание 2.} В зависимости от множества $Q$ показатель $p$ может быть меньше чем $2k$. Например, если $k=2$ и $Q$ достаточно малый квадрат с центром в точке $(1,1),$ то можно доказать, что при $p>3$ интеграл $\theta$ сходится.

\section{Вспомогательные леммы}

Сначала рассматривается следующий несобственный интеграл
$$T_{\infty}(A,b)=\displaystyle\int\limits_{\mathbb{R}^k}{\exp{(iP(x,
A,b)-(x,x))}dx}.$$ Очевидно, что последний интеграл абсолютно и равномерно сходится по параметрам и он явно вычисляется \cite{M}.

\textbf{Лемма 3.}\label{T-imfty}
\emph{Справедливо следующее равенство
$$T_{\infty}(A,b)={(2\pi)}^{\frac{k}{2}}{{(det(I-iA))}^{-\frac{1}{2}}}
{\exp\big({{-\frac{{({{(I-iA)}^{-1}}b,b)}}{4}}}}\big),$$
где квадратный корень определится понимается следующим образом
$${{(det(I-iA))}
^{-\frac{1}{2}}}={(1-i{\lambda}_{1})^{-\frac{1}{2}}}\cdot{(1-i{\lambda}
_{2})^{-\frac{1}{2}}}\cdot...\cdot{(1-i{\lambda}_{k})^{-\frac{1}{2}}}$$
здесь ${{\lambda}_{1}},...,{{\lambda}_{k}}$ собственные значения матрицы $A$.
${z}^{-\frac{1}{2}}$ ветвь многозначной функции, определенной на комплексной плоскости с разрезом по нижней части мнимой оси и ${1}^{-\frac{1}{2}}=1$.}

\textbf{Лемма 3 доказывается} приведением $A$ к диагональному виду. Таким образом вычисление интеграла сводится к одномерному
интегралу и явно вычисляется (более подробно см.\cite{S}).

Очевидно, что выполняются следующие равенства:
$$\left|{\exp\big({{-\frac{{({{(I-iA)}^{-1}}b,b)}}{4}}}}\big)\right|^p=
{\exp({{-\frac{{({{(I+{A}^{2})}^{-1}}b,b)p}}{4}}}})$$ и
$$\displaystyle\int\limits_{R^k}{\exp\big({{-\frac{{({{(I+{A}^{2})}^{-1}}b
,b)p}}{4}}}}\big)db=\frac{{(8\pi)}^{\frac{k}{2}}{
(det(I+{A}^{2}))^{\frac{1}{2}}}}{p^{\frac{k}{2}}}\,.$$

Введем следующее обозначение:
$$\theta_\infty=\displaystyle\int\limits_{{R}^N}{{|T_\infty(A,b)|^{p}}dbda}
,$$ где $$N={\frac{k(k+2)}{2}}.$$
\textbf{Предложение 1.}\label{T-imfty}
\emph{Интеграл \, $\theta_\infty$ сходится при ${p>2k+2}$ и расходится при ${p\leq{2k+2}}$.}

Заметим, что согласно лемме 3, доказательство предложения сводится к исследованию сходимости следующего интеграла
\begin{equation}\label{infty}
\theta_\infty=c(p)\displaystyle\int\limits_{{\mathbb{R}}^{N-k}}{\frac{da}
{{(det(I+{A}^{2}))^{\frac{p-2}{4}}}}},
\end{equation}
где $c(p)$ некоторое положительное число, оно явно вычисляется.

Как известно, определитель является инвариантом ортогональной группы. Поэтому естественно интегрировать сначала по орбитам
ортогональной группы, затем интегрировать по фактор-пространству.

Пусть $M$ множество вещественных симметричных матриц и $SO_{k}$ группа специальных ортогональных матриц. Эта группа естественным образом действует в пространстве $M$, $ g(A)=g^{t}Ag$, где $g{\in{SO_{k}}}$
и $A{\in{M}}$.

Известно, что для любой вещественной симметричной матрицы $A$, существует
$g{\in{G}}$ такое, что $g(A)=diag({{\lambda}_{1}},...,{{\lambda}_{k}
}),$ где $diag({{\lambda}_{1}},...,{{\lambda}_{k}})$ диагональная матрица с диагональными элементами ${{\lambda}_{1}},...,{{\lambda}_{k}}$.
Другими словами, для любой матрицы $A$ существует $g{\in{SO_{k}}}$ такое, что
$ A=g^{t}{\Lambda{g}}$, где
${\Lambda}=diag({{\lambda}_{1}},...,{{\lambda}_{k}})-$некоторая диагональная матрица.

Таким образом, если рассмотреть многообразие ${\mathbb{R}^{k}{\times{SO_{k}}}}$, то естественно определяется гладкое сюръективное отображение
$${\Phi}:{\mathbb{R}^{k}{\times{SO_{k}}}}{\mapsto}M$$ определенное по формуле ${\Phi({\Lambda},g)}=g^{t}{\Lambda{g}}$.

Пусть $da=da_{11}{\wedge}da_{12}{\wedge}...{\wedge}da_{kk}$ естественная форма объема в пространстве $M$. Мы можем определить образ этой формы при отображении
${\Phi}$, обозначаемый через
${\Phi}^{\ast}da{\in}{{\wedge}^{N-k}}({\mathbb{R}^{k}{\times{SO_{k}}}})$.

\textbf{Лемма 4.}\label{Phi}
\emph{Справедливо следующее равенство
$${\Phi}^{\ast}da=\prod\limits_{{1{\leq{l}}<{m}{\leq{k}}}}
({\lambda}_{m}-{\lambda}_{l})d{\lambda}_{1}{\wedge}...{\wedge}
d{\lambda}_{k}{\wedge}{\omega},$$ где ${\omega}-$форма объема на ортогональной группе $SO_{k}$.}

\textbf{Лемма 4 доказывается} с использованием нулевого множества якобиана отображения ${\Phi}$. Отметим, что справедливо равенство
$\prod\limits_{{1{\leq{l}}<{m}{\leq{k}}}}({\lambda}_{m}-{\lambda}_{l})
^{2}={\rho}_{A}({\lambda})$, где
${\rho}_{A}({\lambda})-$ха\-рак\-те\-рис\-ти\-чес\-кий многочлен матрицы $A$.

Согласно лемме 3, интеграл (\ref{infty}) записывается в виде
$${\displaystyle\int\limits_{\mathbb{R}^{N-k}}{\frac{da}{{(det(I+{A}^{2}))
^{\frac{p-2}{4}}}}}}=\displaystyle\int\limits_{\mathbb{R}^{k}}{\frac{\prod
\limits_{{1{\leq{l}}<{m}{\leq{k}}}}{|{\lambda}_{m}-{\lambda}_{l}|}}{\prod
\limits_{{1{\leq{l}}{\leq{k}}}}(1+{\lambda_{l}}^{2})^{\frac{p-2}{4}}}{d{
\lambda}_{1}{\wedge}...{\wedge}d{\lambda}_{k}}}{\displaystyle\int\limits_{SO_{k}}
{\omega}}$$ Из последнего равенства следует, что сходимость интеграла (\ref{infty}) сводится к исследованию сходимости интеграла
$${\displaystyle\int\limits_{\mathbb{R}^{k}}{\frac{\prod\limits_{{1{\leq{l}}
<{m}{\leq{k}}}}{|{\lambda}_{m}-{\lambda}_{l}|}}{\prod\limits_{{1{\leq{l}}{
\leq{k}}}}(1+{\lambda_{l}}^{2})^{\frac{p-2}{4}}}{d{\lambda}_{1}{\wedge}...
{\wedge}d{\lambda}_{k}}}}.$$ Легко видеть, что этот интеграл сходится при
${p>2k+2}$ и расходится при ${p\leq{2k+2}}$. Что и доказывает предложение 1.

\textbf{Доказательство теоремы 2.}
При доказательстве используется классическое неравенство Юнга.
Если ${f{\in{L_{p}(\mathbb{R}^{k})}}}$ и
${g{\in{L_{r}(\mathbb{R}^{k})}}}$ произвольные функции, то справедливо следующее неравенство
$$\|f{\ast}g\|_{L_{q}}{\leq}\|{f}\|_{L_{p}}\|{g}\|_{L_{r}},$$
где  $f{\ast}g$ свертка функции $f$ и $g$, причем
постоянные $1{\leq}{p,q,r}{\leq}{\infty}$ связаны соотношением
$${{\frac{1}{q}}+1={\frac{1}{p}}+{\frac{1}{r}}}.$$

Пусть $K$ компактный многогранник в $\mathbb{R}^{k}$ и
$$h(b)={\displaystyle\int\limits_{\mathbb{R}^{k}}
{e^{{\mid{x}{\mid}}^{2}}}{{\chi}_{K}}(x){e^{-2\pi i(b,x)}}dx}.$$

\textbf{Лемма 5.}\label{h}
\emph{Для любого положительного числа ${\epsilon}$, имеет место включение
$h{\in{L_{1+{\epsilon}}({\mathbb{R}^{k}})}}$.}

\textbf{Доказательство леммы 5.} Заметим, что для любого $\varepsilon>0$ $\hat{\chi}_{K}\in\mathbb{L}_{1+\varepsilon}(\mathbb{R}^{k})$ (например, см.\cite{Lebedev}).
тогда утверждение леммы 3 легко следует из неравенства Юнга.

Теперь вернемся к доказательству теоремы 2. Согласно тождеству Планшареля имеем:
$$T(A)=\int_{K}e^{i(Ax,x)}dx=\int_{\mathbb{R}^{k}}e^{i(Ax,x)}\chi_{K}(x)dx=
\int_{\mathbb{R}^{k}}e^{i(Ax,x)-|x|^{2}}e^{|x|^{2}}\chi_{K}(x)dx=\int_{\mathbb{R}^{k}}\widehat{f}(A,b)\overline{g}(b)db,$$
где
$\widehat{f}(A,b)=\int_{\mathbb{R}^{k}}e^{i(Ax,x)-|x|^{2}-2\pi i(x,b)}dx$ и $\widehat{g}(b)=\int_{\mathbb{R}^{k}}e^{|x|^{2}}e^{-2\pi i(x,b)}dx$.

 Пусть $q>1$ фиксированное число. Тогда, применяя неравенство Гёльдера, имеем:
$$|T(A)|\leq\|\widehat{f}(A,\cdot)\|_{\mathbb{L}^{q^{'}}(\mathbb{R}^{k})}\|g\|_{\mathbb{L}^{q}(\mathbb{R}^{k})},$$
где $\frac{1}{q}+\frac{1}{q^{'}}=1.$

Согласно лемме 3 имеем:
$$|T(A)|\leq\frac{c_{q}}{(det(I+A^{2}))^{\frac{p}{4}-\frac{1}{2q^{'}}}}.$$
Таким образом, если $p>2k,$ то мы можем выбрать $q^{'}>1$ так, что $\frac{p}{4}-\frac{1}{2q^{'}}>\frac{k}{2}.$

Отсюда следует, что если $\frac{p}{4}-\frac{1}{2q^{'}}>\frac{k}{2},$ то $T\in\mathbb{L}^{p}(\mathbb{R}^{k}).$

Осталось доказать точность результата. Рассмотрим следующее подмножество
${\Omega^{+}{(a_{11})}}$ пространства $\mathbb{R}^{N-1},$ где $N=\frac{k(k+1)}{2}.$
$$a_{11}>0, {|a_{12}|+|a_{13}|+...+|a_{1k}|}<c_{1}a_{11}, \quad
\left|{a_{lj}}-{\frac{{a}_{1l}{a_{1j}}}{a_{11}}}\right|{\leq{c_{2}}}, a_{1l}<0$$ где
$l\leq j=\overline{2,n}$,$l=2,...,n$ и $c_{1}$, $c_{2}$ достаточно малые фиксированные положительные числа.

Согласно лемме 1 существуют положительные числа $c_{1}$ и $c_{2}$ такие, что для меры Лебега множества ${\Omega^{+}{(a_{11})}}$ справедливо следующее равенство:
$${\mu{({\Omega^{+}{(a_{11})}})}} =c\cdot{a^{k-1}_{11}}.$$

\textbf{Лемма 6.}\label{T(A,b)}
\emph{Существуют положительное число $L$ такое, что для любого ${a_{11}}>L$ и
$A{\in{\Omega^{+}{(a_{11})}}}$ для интеграла $T(A)$ справедливо следующее асимптотическое равенство
$$T(A)={\frac{c(A)}{a^{\frac{1}{2}}_{11}}+O\left({\frac{1}{a_{11}}}\right)}$$
 при ${a_{11}}{\rightarrow}{+{\infty}},$
причем существует положительное число ${\delta}$ такое, что для любого
$(A,b){\in{{\Omega^{+}{(a_{11})}}}}$ выполняется неравенство
$$|c(A)|>{\delta}>0.$$}

\textbf{Лемма 6 доказывается} обычным методом стационарной фазы. Заметим, что если $\delta_{2}>\frac{1}{2}$ и $\delta_{1}<0$ то справедливо следующее соотношение
$$\left|\int_{\delta_{1}\sqrt{\lambda}}^{\delta_{2}\sqrt{\lambda}}\cos y^{2}dy\right|=c(\delta_{1},\delta_{2},\lambda)$$ причем, существуют $\lambda_{0}$, $\varepsilon>0$ такие, что выполняется неравенство $c(\delta_{1},\delta_{2},\lambda)\geq\varepsilon>0$ при всех $\lambda\geq\lambda_{0}.$

Отметим, что для достаточно малых $c_{1}$, $c_{2}$ при $A\in{\Omega^{+}{(a_{11})}}$ и для достаточных больших $L$ фаза имеет осцилляции только в направлении $x_{1}$ по этому, при фиксированных значениях  ${x_{2},...,x_{n}{\in{[0,1]}}}$ невырожденная критическая точка ${x_{1}}(A,b,x_{2},...,x_{n})$ лежит внутри $(0,1).$

Наконец, для интеграла ${\theta}$, имеем оценку снизу
$${\theta}{\geq}{\displaystyle\int\limits_{L}^{\infty}}
{\displaystyle\int\limits_{{\Omega{(a_{11})}}}{{|T(A)|^{p}}da}}
{\geq}{\delta}{c}{\displaystyle\int\limits_{L}^{\infty}}{a^{k-
{\frac{p}{2}-1}}_{11}}da_{11}.$$

Таким образом при ${p\leq{2k}}$ последний интеграл расходится. Основная теорема 2 доказана.

\section{Двумерный случай}
Отметим, что в однородном случае результаты \cite{LeeBak} не применимы. При доказательстве теоремы 2 существенно используется свойство $\widehat{\chi}_{Q}\in\mathbb{L}_{1+0}(\mathbb{R}^{k}).$

В работе В.В.Лебедева приведен пример области $\partial D\in C^{1,\omega}$, где $\omega$ модуль непрерывности градиента $\varphi$, определяющей $\partial D.$ Поэтому, мы можем считать, что $D$ компактная область.

Справедлива следующая
\begin{theorem}
  Пусть $D$ компактная область такая, что $\widehat{\chi}_{D}\in\mathbb{L}_{q}(\mathbb{R}^{2})$ и
  $T(A)=\int_{D}e^{i(Ax,x)}dx.$ Тогда $T\in\mathbb{L}^{p}(\mathbb{R}^{3})$ при $p>6-\frac{2}{q}.$ Более того, если
  $\widehat{\chi}_{D}\in\mathbb{L}_{1+0}(\mathbb{R}^{2}),$ то, при любом $p>4$, справедливо включение $T\in\mathbb{L}^{p}(\mathbb{R}^{3}).$
\end{theorem}

\textbf{Замечание 3.} Из результатов В.В.Лебедева \cite{Lebedev} следует, что существует множество
$D$ отличное от многоугольника, такое, что $\widehat{\chi}_{Q}\in\mathbb{L}_{1+0}(\mathbb{R}^{2}).$

\textbf{Следствие.} Если $D\subset\mathbb{R}^{2}$ компактное множество, такое, что $\partial D\subset C^{1}$, то при $p>4,5$ справедливо соотношение $T\in\mathbb{L}^{p}(\mathbb{R}^{3}).$

%библиография по ГОСТу

\end{document}